\documentclass[twoside,twocolumn]{article}
\usepackage{latexsym}
\usepackage{graphicx}
\makeatletter
\def\section{\@startsection {section}{1}{\z@}{-2.0ex plus
    -0.5ex minus -.2ex}{1.5ex plus 0.3ex minus .2ex}{\large\bf\raggedright}}
\makeatother

\def\,{\mskip 3mu} \def\>{\mskip 4mu plus 2mu minus 4mu} \def\;{\mskip 5mu plus 5mu} \def\!{\mskip-3mu}
\def\dispmuskip{\thinmuskip= 3mu plus 0mu minus 2mu \medmuskip=  4mu plus 2mu minus 2mu \thickmuskip=5mu plus 5mu minus 2mu}
\def\textmuskip{\thinmuskip= 0mu                    \medmuskip=  1mu plus 1mu minus 1mu \thickmuskip=2mu plus 3mu minus 1mu}
\textmuskip
\def\beq{\dispmuskip\begin{equation}}    \def\eeq{\end{equation}\textmuskip}
\def\beqn{\dispmuskip\begin{displaymath}}\def\eeqn{\end{displaymath}\textmuskip}
\def\bqa{\dispmuskip\begin{eqnarray}}    \def\eqa{\end{eqnarray}\textmuskip}
\def\bqan{\dispmuskip\begin{eqnarray*}}  \def\eqan{\end{eqnarray*}\textmuskip}

\def\paradot#1{\vspace{1.3ex plus 0.5ex minus 0.5ex}\noindent{\bf{#1.}}}
\def\paranodot#1{\vspace{1.3ex plus 0.5ex minus 0.5ex}\noindent{\bf{#1}}}

\def\aidx#1{}

\def\eps{\varepsilon}
\def\epstr{\epsilon}                    
\def\nq{\hspace{-1em}}

\def\odt{{\textstyle{1\over 2}}}
\def\SetR{I\!\!R}
\def\SetN{I\!\!N}
\def\SetB{I\!\!B}
\def\SetZ{Z\!\!\!Z}
\def\qmbox#1{{\quad\mbox{#1}\quad}}
\def\e{{\rm e}}                        
\def\B{\{0,1\}}
\def\v{\vec}
\def\es{\phi}                          
\def\l{\ell}
\def\lb{{\log_2}}

\begin{document}

\title{\normalsize\sc Technical Report \hfill IDSIA-24-04
\vskip 2mm\bf\Large\hrule height5pt \vskip 6mm
Fast Non-Parametric Bayesian Inference on Infinite Trees
\vskip 6mm \hrule height2pt \vskip 5mm}
\author{{\bf Marcus Hutter}\\[3mm]
\normalsize IDSIA, Galleria 2, CH-6928\ Manno-Lugano, Switzerland \\
\normalsize marcus@idsia.ch, \ http://www.idsia.ch/$^{_{_\sim}}\!$marcus }
\date{21 May 2004}
\maketitle

\begin{abstract}
Given i.i.d.\ data from an unknown distribution, we consider the
problem of predicting future items. An adaptive way to estimate
the probability density is to recursively subdivide the domain to
an appropriate data-dependent granularity. A Bayesian would assign
a data-independent prior probability to ``subdivide'', which leads
to a prior over infinite(ly many) trees. We derive an exact, fast,
and simple inference algorithm for such a prior, for the data
evidence, the predictive distribution, the effective model
dimension, and other quantities.
\end{abstract}

\section{INTRODUCTION}\label{secInt}

\paradot{Inference}
We consider the problem of inference from i.i.d.\ data $D$, in
particular of the unknown distribution $q$ the data is sampled
from. In case of a continuous domain this means inferring a
probability density from data. Without structural assumption on
$q$, this is hard to impossible, since a finite amount of data is
never sufficient to uniquely select a density (model) from an
infinite-dimensional space of densities (model class).

\paradot{Methods}
In parametric estimation one assumes that $q$ belongs to a
finite-dimensional family. The two-dimensional family of Gaussians
characterized by mean and variance is prototypical. The maximum
likelihood (ML) estimate of $q$ is the distribution that
maximizes the data likelihood. Maximum likelihood overfits if the
family is too large and especially if it is infinite-dimensional.
A remedy is to penalize complex distributions by assigning a prior
(2nd order) probability to the densities $q$. Maximizing the model
posterior (MAP), which is proportional to likelihood times the prior,
prevents overfitting.
Bayesians keep the complete posterior for inference. Typically,
summaries like the mean and variance of the posterior are
reported.

\paranodot{How to choose the prior?}
In finite or small compact low-dimensional spaces a uniform prior
often works (MAP reduces to ML). In the non-parametric case one
typically devises a hierarchy of finite-dimensional model classes
of increasing dimension. Selecting the dimension with maximal
posterior often works well due to the Bayes factor phenomenon
\cite{Good:83,Jaynes:03,MacKay:03}: In case the true model is
low-dimensional, higher-dimensional (complex) model classes are
automatically penalized, since they contain fewer ``good'' models.
Full Bayesians would assign a prior probability (e.g.\ ${1\over
d^2}$) to dimension $d$  and mix over dimension.

\paradot{Interval Bins}
The probably simplest and oldest model for an interval domain is
to divide the interval (uniformly) into bins, assume a constant
distribution within each bin, and take a frequency estimate for
the probability in each bin, or a Dirichlet posterior if you are a
Bayesian. There are heuristics for choosing the number of bins as
a function of the data size. The simplicity and easy computability
of the bin model is very appealing to practitioners. Drawbacks are
that distributions are discontinuous, its restriction to one
dimension (or at most low dimension: curse of dimensionality), the
uniform (or more generally fixed) discretization, and the
heuristic choice of the number of bins. We present a full Bayesian
solution to these problems, except for the non-continuity problem.
Polya trees \cite{Lavine:94} inspired our model.

\paradot{More advanced model classes}
There are plenty of alternative Bayesian models that
overcome some or all of the limitations. Examples are %
continuous Dirichlet process (mixtures) \cite{Ferguson:73}, %
Bernstein polynomials \cite{Petrone:02}, %
Bayesian field theory \cite{Lemm:03}, %
Bayesian kernel density estimation %
or other mixture models \cite{Escobar:95},
or universal priors \cite{Hutter:04uaibook}, %
but analytical solutions are infeasible. %
Markov Chain Monte Carlo sampling %
or Expectation Maximization algorithms \cite{Dempster:77} %
or variational methods %
can often be used to obtain approximate numerical solutions, but
computation time and global convergence remain critical issues.
Practitioners usually use (with success) efficient MAP or M(D)L or
heuristic methods, e.g.\ kernel density estimation \cite{Gray:03},
but note that MAP or MDL {\em can} fail, while Bayes works
\cite{Hutter:04mdl2p}.

\paradot{Our tree mixture model}
The idea of the model class discussed in this paper is very
simple:
With equal probability, we chose $q$ either uniform or
split the domain in two parts (of equal volume), and assign a
prior to each part, recursively, i.e.\ in each part again either
uniform or split. For finitely many splits, $q$ is a
piecewise constant function, for infinitely many splits it is
virtually {\em any} distribution. While the prior over $q$ is
neutral about uniform versus split, we will see that the posterior
favors a split if and only if the data clearly indicates
non-uniformity. The method is a full Bayesian non-heuristic
tree approach to adaptive binning for which we present a very
simple and fast algorithm for computing all(?) quantities of
interest.

\paradot{Contents}
In Section \ref{secTMM} we introduce our model and compare it to
Polya trees. We also discuss some example domains, like intervals,
strings, volumes, and classification tasks.
In Section \ref{secQI} we present recursions for various
quantities of interest, including the data evidence, the
predictive distribution, the effective model dimension, the tree
size and height, and cell volume. We discuss the qualitative
behavior and state convergence of the posterior for finite trees.
The proper case of infinite trees is discussed in Section
\ref{secIT}, where we analytically solve the infinite recursion at
the data separation level.
Section \ref{secAlg} collects everything together and presents the
algorithm.
We also numerically illustrate the behavior
of our model on one example distribution.
Section \ref{secDisc} contains a brief summary, conclusions, and
outlook, including natural generalizations of our model.
See \cite{Hutter:04btcode} for derivations, proofs, program code,
extensions, and more details.

\section{THE TREE MIXTURE MODEL}\label{secTMM}

\paradot{Setup and basic quantities of interest}
We are given i.i.d.\ data $D=(x^1,...,x^n)\in\Gamma^n$ of size $n$
from domain $\Gamma$, e.g.\ $\Gamma\subseteq\SetR^d$ sampled from
some unknown probability density $q:\Gamma\to\SetR$. Standard
inference problems are to estimate $q$ from $D$ or to predict the
next data item $x^{n+1}\in\Gamma$. By definition, the (objective
or aleatoric) data likelihood density under model $q$ is
$p(D|q)\equiv q(x_1)\cdot...\cdot q(x_n)$. Note that we consider
sorted data, which avoids annoying multinomial coefficients.
Otherwise this has no consequences. Results are independent of the
order and depend on the counts only, as they should. A Bayesian
assumes a (belief or $2^{nd}$-order or epistemic or subjective)
prior $p(q)$ over models $q$ in some model class $Q$. The data
evidence is $p(D)=\int_Q p(D|q)p(q) dq$. Having the evidence,
Bayes' famous rule allows to compute the (belief or $2^{nd}$-order
or epistemic or subjective) posterior $p(q|D)={p(D|q)p(q)/p(D)}$
of $q$. The predictive or posterior distribution of $x$ is
$p(x|D)={p(D,x)/p(D)}$, i.e.\ the conditional probability that the
next data item is $x=x^{n+1}$, given $D$, follows from the
evidences of $D$ and $(D,x)$. Since the posterior of $q$ is a
complex object, we need summaries like the expected
$q$-probability of $x$ and (co)variances. Fortunately they can
also be reduced to computation of evidences: ${E[q(x)|D]} := {\int
q(x)p(q|D) dq} = {p(x|D)}$. In the last equality we used the
formulas for the posterior, the likelihood, the evidence, and the
predictive distribution, in this order. Similarly for the
covariance. We derive and discuss further summaries of $q$ for our
particular tree model, like the model complexity or effective
dimension, and the tree height or cell size, later.

\paradot{Hierarchical tree partitioning}
Up to now everything has been fairly general. We now introduce the
tree representation of domain $\Gamma$. We partition $\Gamma$ into
$\Gamma_0$ and $\Gamma_1$, i.e.\ $\Gamma=\Gamma_0\cup\Gamma_1$ and
$\Gamma_0\cap\Gamma_1=\es$. Recursively we (sub)partition
$\Gamma_z=\Gamma_{z0}\dot\cup\Gamma_{z1}$ for $z\in\SetB_0^m$,
where $\SetB_k^m=\bigcup_{i=k}^m\{0,1\}^i$ is the set of all
binary strings of length between $k$ and $m$, and
$\Gamma_\epstr=\Gamma$, where $\epstr=\{0,1\}^0$ is the empty
string. We are interested in an infinite recursion, but for
convenience we assume a finite tree height $m<\infty$ and
consider $m\to\infty$ later. Also let $l:=\l(z)$ be the length of
string $z=z_1...z_l=:z_{1:l}$, and $|\Gamma_z|$ the
volume or length or cardinality of $\Gamma_z$.

\paradot{Example spaces}
{\it Intervals:} Assume $\Gamma=[0,1)$ is the unit interval, recursively
bisected into intervals $\Gamma_z=[0.z, 0.z+2^{-l})$ of length
$|\Gamma_z|=2^{-l}$, where $0.z$ is the real number in $[0,1)$
with binary expansion $z_1...z_l$.

{\it Strings:} Assume $\Gamma_z=\{zy:y\in\B^{m-l}\}$ is the set of
strings of length $m$ starting with $z$. Then $\Gamma=\B^m$ and
$|\Gamma_z|=2^{m-l}$. For $m=\infty$ this set is continuous, for
$m<\infty$ finite.

{\it Trees:} Let $\Gamma$ be a complete binary tree of height $m$ and
$\Gamma_{z0}$ ($\Gamma_{z1}$) be the left (right) subtree of
$\Gamma_z$. If $|\Gamma_z|$ is defined as one more than the number
of nodes in $\Gamma_z$, then $|\Gamma_z|=2^{m+1-l}$.

{\it Volumes:} Consider $\Gamma\subset\SetR^d$, e.g.\ the
hypercube $\Gamma=[0,1)^d$. We recursively halve $\Gamma_z$ with a
hyperplane orthogonal to dimension $(l$ mod $d)+1$, i.e.\ we sweep
through all orthogonal directions. $|\Gamma_z|=2^{-l}|\Gamma|$.

{\it Compactification:} We can compactify
$\Gamma\subseteq(1,\infty]$ (this includes
$\Gamma=\SetN\setminus\{1\}$) to the unit interval
$\Gamma':=\{{1\over x}:x\in\Gamma\}\subseteq[0,1)$, and similarly
$\Gamma\subseteq\SetR$ (this includes $\Gamma=\SetZ$) to
$\Gamma':=\{x\in[0,1):{2x-1\over x(1-x)}\in\Gamma\}$. All
reasonable spaces can be reduced to one of the spaces described above.

{\it Classification:} Consider an observation $o\in\Gamma'$ (e.g.\
email) that is classified as $c\in\{0,1\}$ (e.g.\ good versus
spam), where $\Gamma'$ could be one of the spaces above (e.g.\
$o$ is a sequence of binary features in decreasing order of
importance). Then $x:=(o,c)\in\Gamma:=\Gamma'\times\{0,1\}$ and
$\Gamma_{0z}=\Gamma'_z\times\{0\}$ and
$\Gamma_{1z}=\Gamma'_z\times\{1\}$. Given $D$ (e.g.\
pre-classified emails), a new observation $o$ is classified as $c$
with probability $p(c|D,o)\propto p(D,x)$. Similar for more than
two classes.

In all these examples (we have chosen)
$|\Gamma_{z0}|=|\Gamma_{z1}|=\odt|\Gamma_z|$ $\forall
z\in\SetB_0^{m-1}$, and this is the only property we need and
henceforth assume. W.l.g.\ we assume/\linebreak[1]define/\linebreak[1]rescale $|\Gamma|=1$.
Generalizations to non-binary and non-symmetric partitions are
straightforward and briefly discussed at the end.

\paradot{Identification}
We assume that $\{\Gamma_z:z\in\SetB_0^m\}$ are (basis) events
that generate our $\sigma$-algebra. For every $x\in\Gamma$ let
$x'$ be the string of length $\l(x')=m$ such that $x\in\Gamma_{x'}$.
We assume that distributions $q$ are $\sigma$-measurable, i.e. to
be constant on $\Gamma_{x'}$ $\forall x'\in\SetB^m$. For
$m=\infty$ this assumption is vacuous; we get {\em all} Borel
measures. Hence, we can identify the continuous sample space
$\Gamma$ with the (for $m<\infty$ discrete) space $\SetB^m$ of
binary sequences of length $m$, i.e.\ in a sense all example
spaces are isomorphic.
While we have the volume model in mind for real-world
applications, the string model will be convenient for mathematical
notation, the tree metaphor will be convenient in discussion, and
the interval model will be easiest to implement and to present
graphically.

\paradot{Notation}
As described above, $\Gamma$ may also be a tree. This
interpretation suggests the following scheme for defining the
probability of $q$ on the leaves $x'$. The probability of the left
child node $z0$, given we are in the parent node $z$, is
$P[\Gamma_{z0}|\Gamma_z,q]$, so we have
\beqn
   p(x|\Gamma_z,q) = p(x|\Gamma_{z0},q)\!\cdot\! P[\Gamma_{z0}|\Gamma_z,q]
   \qmbox{if} x\in\Gamma_{z0}
\eeqn
and similarly for the right child. In the following we often have
to consider distributions conditioned to and in the subtree
$\Gamma_z$, so the following notation will turn out
convenient\vspace{-1ex}
\beq\label{eqNotation}
  q_{z0} := P[\Gamma_{z0}|\Gamma_z,q], \quad
  p_z(x|...) := 2^{-l}p(x|\Gamma_z...)\vspace{-1ex}
\eeq
\beqn
  \Rightarrow p_z(x|q) = 2q_{zx_{l+1}} p_{zx_{l+1}}(x|q) =...
  = \!\!\prod_{i=l+1}^m\! 2q_{x_{1:i}}
  \;\mbox{if}\; x\in\Gamma_z
\eeqn
where we have used that $p(x|\Gamma_{x'},q)=|\Gamma_{x'}|^{-1}=2^m$ is
uniform. Note that $q_{z0}+q_{z1}=1$. Finally, let
$\v q_{z*}:=(q_{zy}:y\in\SetB_1^{m-l})$
be the ($2^{m-l+1}-2$)-dimensional {\em vector} or {\em ordered
set} or {\em tree} of all {\em reals} $q_{zy}\in[0,1]$ in subtree
$\Gamma_z$. Note that $q_z\not\in\v q_{z*}$. The {\em
(non)density} $q_z(x):= p_z(x|q)$ depends on all and only these
$q_{zy}$. For $z\neq\epstr$, $q_z()$ and $p_z()$ are only
proportional to a density due to the factor $2^{-l}$, which has
been introduced to make $p_{x'}(x|...)\equiv 1$. (They are
densities w.r.t.\ $2^l\lambda_{|\Gamma_z}$, where $\lambda$ is the
Lebesgue measure.) We have to keep this in mind in our
derivations, but can ignore this widely in our discussion.

\paradot{Polya trees}
In the Polya tree model one assumes that the $q_{z0}\equiv
1-q_{z1}$ are independent and Beta($\cdot,\cdot$) distributed,
which defines the prior over $q$. Polya trees form a conjugate
prior class, since the posterior is also a Polya tree, with
empirical counts added to the Beta parameters. If the same Beta is
chosen in each node, the posterior of $x$ is pathological for
$m\to\infty$: The distribution is everywhere discontinuous with
probability 1. A cure is to increase the Beta parameters with $l$,
e.g.\ quadratically, but this results in ``underfitting'' for
large sample sizes, since Beta(large,large) is too informative and
strongly favors $q_{z0}$ near $\odt$. It also violates scale
invariance, which should hold in the absence of prior knowledge.
That is, the p(oste)rior in $\Gamma_0=[0,\odt)$ should be the same
as for $\Gamma=[0,1)$ (after rescaling all $x\leadsto x/2$ in
$D$).

\paradot{The new tree mixture model}
The prior $P[q]$ follows from specifying a prior over $\v q_*$, since
$q(x)\propto q_{x_1}\cdot...\cdot q_{x_{1:m}}$ by (\ref{eqNotation}). The
distribution in each subset $\Gamma_z\subseteq\Gamma$ shall be
either $u$niform or non-uniform. A necessary (but not sufficient)
condition for uniformity is $q_{z0}=q_{z1}=\odt$.
\beq\label{tmmu1}
  p^u(q_{z0},q_{z1}) \;:=\; \delta(q_{z0}-\odt)\delta(q_{z1}-\odt),
\eeq
where $\delta()$ is the Dirac delta. To get uniformity on
$\Gamma_z$ we have to recurse the tree down in this way.
\beq\label{tmmur}
  p^u(\v q_{z*}) \;:=\; p^u(q_{z0},q_{z1})p^u(\v q_{z0*})p^u(\v q_{z1*})
\eeq
with the natural recursion termination $p^u(\v q_{z*})=1$ when
$\l(z)=m$, since then $\v q_{z*}=\es$.
For a non-uniform distribution on $\Gamma_z$ we allow
any probability split $q(\Gamma_z)=q(\Gamma_{z0})+q(\Gamma_{z0})$,
or equivalently $1=q_{z0}+q_{z1}$. We assume a uniform prior on the
$s$plit, i.e.
\beq\label{tmms1}
   p^s(q_{z0},q_{z1}) \;:=\; \delta(q_{z0}+q_{z1}-1)
\eeq
We now recurse down the tree
\beq\label{tmmsr}
   p^s(\v q_{z*}) \;:=\; p^s(q_{z0},q_{z1})p(\v q_{z0*})p(\v q_{z1*})
\eeq
again with the natural recursion termination $p(\v q_{z*})=p(\es)=1$
when $\l(z)=m$. Finally we have to mix the uniform with the
non-uniform case.
\beq\label{tmmmix}
  p(\v q_{z*}) \;:=\; p(u)p^u(\v q_{z*})+p(s)p^s(\v q_{z*})
\eeq
We choose a 50/50 mixture $p(u)=p(s)=\odt$. This completes the
specification of the prior $P[q]=p(\v q_*)$.

For example, if the first bit in $x$ is a class label
and the remaining are binary features in decreasing order of
importance, then given class and features $z=x_{1:l}$, further features
$x_{l+1:m}$ could be relevant for classification ($q_z(x)$ is
non-uniform) or irrelevant ($q_z(x)$ is uniform).

\paradot{Comparison to the Polya tree}
Note the important difference in the recursions (\ref{tmmur}) and
(\ref{tmmsr}). Once we decided on a uniform distribution
(\ref{tmmu1}) we have to equally split probabilities down the recursion to the
end, i.e.\ we recurse in (\ref{tmmur}) with $p^u$, rather than the
mixture $p$ (this actually allows to solve the recursion). On the
other hand if we decided on a non-uniform split
(\ref{tmms1}), the left and right partition each itself may be
uniform or not, i.e.\ we recurse in (\ref{tmmsr}) with the mixture
$p$, rather than $p^s$. Inserting (\ref{tmms1}) in (\ref{tmmsr})
in (\ref{tmmmix}) and recursively (\ref{tmmu1}) in (\ref{tmmur})
in (\ref{tmmmix}) we get
\beq\label{EqPqRec}
  p(\v q_{z*}) = \odt\nq\smash{\prod_{y\in\SetB_1^{m-l}}}\nq \delta(q_{zy}\!-\!\odt)
  + \odt\delta(q_{z0}\!+\!q_{z1}\!-\!1)p(\v q_{z0*})p(\v q_{z1*})
\eeq
Choosing $p(u)=0$ would lead to the Polya tree model (and its
problems) with $q_{z0}\sim$ Beta(1,1). For our choice
($p(u)=\odt$), but with $p$ instead of $p^u$ on the r.h.s.\ of
(\ref{tmmur}) we would get a quasi-Polya model (same problems) with
$q_{z0}\sim\odt[\mbox{Beta}(\infty,\infty)+\mbox{Beta}(1,1)]$.

For $m\to\infty$, our model is scale invariant {\em and} leads to
continuous distributions for $n\to\infty$, unlike the Polya tree
model. We also don't have to tune Beta parameters; the model tunes
itself by suitably assigning high/low posterior probability to
subdividing cells. While Polya trees form a natural conjugate
prior class, our prior does not directly, but can be generalized
to do so \cite{Hutter:04btcode}. The computational complexity
for the quantities of interest will be the same (essentially
$O(n)$), i.e.\ as good as it could be.

\paradot{Formal and effective dimension}
Formally our model is $2\cdot(2^m-1)$-dimensional, but the
effective dimension can by much smaller, since $\v q_*$ is forced
with a non-zero probability to a much smaller polytope, for
instance with probability $\odt$ to the zero-dimensional globally
uniform distribution. We will compute the effective p(oste)rior
dimension.

\section{QUANTITIES OF INTEREST}\label{secQI}

\paradot{The evidence recursion}
At the end of Section \ref{secTMM} we defined our tree mixture
model. The next step is to compute the standard quantities of
interest defined at the beginning of Section \ref{secTMM}. The
evidence $p(D)$ is key, the other quantities (posterior,
predictive distribution, expected $q(x)$ and its variance) follow
then immediately.
Let $D_z:=\{x\in D : x\in\Gamma_z\}$ be the $n_z:=|D_z|$ data
points that lie in subtree $\Gamma_z$. We compute $p_z(D_z)$
recursively for all $z\in\SetB_0^{m-1}$, which gives
$p(D)=p_\epstr(D_\epstr)$. Inserting (\ref{eqNotation}) and
(\ref{EqPqRec}) into
\beq\label{eqEvRecProof}
  p_z(D_z) = \int p_z(D_z|\v q_{z*}) p(\v q_{z*}) d{\v q}_{z*}
\eeq
one can derive the following recursion \cite{Hutter:04btcode}:
\bqa\label{eqEvDens}\label{eqWeights}\label{eqEvRec}
  p_z(D_z) &=& \odt\big[1+{p_{z0}(D_{z0})p_{z1}(D_{z1})\over w({n_{z0},n_{z1}})}\big]
\\ \nonumber
  \nq\nq w({n_{z0},n_{z1}}) &:=& 2^{-n_z}{(n_z\!+\!1)!\over n_{z0}!n_{z1}!}
  \;=:\; w_{n_z}(\Delta_z)
\\ \nonumber
  \qquad n_z &=& n_{z0}+n_{z1}, \quad
  \Delta_z \;:=\; {n_{z0}\over n_z}-\odt
\eqa
The recursion terminates with $p_z(D_z)=1$ when $\l(z)=m$. Recall
(\ref{eqNotation}) if you insist on a formal proof: For $\l(z)=m$
and $x\in\Gamma_z$ we have $\Gamma_{x'}=\Gamma_z$ $\Rightarrow$
$p_z(x|q)=1$ $\Rightarrow$ $p_z(D_z|q)=1$ $\Rightarrow$
$p_z(D_z)=1$.

Interpretation of (\ref{eqEvDens}): With probability $\odt$, the
evidence is uniform in $\Gamma_z$. Otherwise data $D_z$ is split
into two partitions of size $n_{z0}$ and $n_{z1}=n_z-n_{z0}$. First,
choose $n_{z0}$ uniformly in $\{0,...,n_z\}$. Second, given $n_z$,
choose uniformly among the $({n_z\atop n_{z0}})$ possibilities of
selecting $n_{z0}$ out of $n_z$ data points for $\Gamma_{z0}$ (the
remaining $n_{z1}$ are then in $\Gamma_{z1}$). Third, distribute
$D_{z0}$ according to $p_{z0}(D_{z0})$ and $D_{z1}$ according to
$p_{z1}(D_{z1})$. Then, the evidence in case of a split is the
second term in (\ref{eqEvRec}). The factor $2^{n_z}$ is due
to our normalization convention (\ref{eqNotation}). This also verifies
that the r.h.s.\ yields the l.h.s.\ if integrated over all $D_z$,
as it should be.

\paradot{\boldmath Discussing the weight}
The relative probability of splitting (second term on r.h.s.\ of
(\ref{eqEvDens})) to the uniform case (first term in r.h.s.\ of
(\ref{eqEvDens})) is controlled by the weight $w$. Large (small)
weight indicates a (non) uniform distribution, provided $p_{z0}$
and $p_{z1}$ are $O(1)$. Balance $\Delta_z\approx 0$ ($\not\approx 0$)
indicates a (non) symmetric partitioning of the data among the left
and right branch of $\Gamma_z$. Asymptotically for large $n_z$
(keeping $\Delta_z$ fixed), we have%
\beqn
  w_{n_z}(\Delta_z) \sim {\textstyle\sqrt{2n_z\over\pi}}\;\e^{-2n_z\Delta_z^2}
\eeqn
Assume that data $D$ is sampled from the true distribution $\dot
q$. The probability of the left branch $\Gamma_{z0}$ of $\Gamma_z$
is $\dot q_{z0}\equiv P[\Gamma_{z0}|\Gamma_z,\dot q]=
2^l\dot q_z(\Gamma_{z0})$. The relative frequencies ${n_{z0}\over
n_z}$ asymptotically converge to $\dot q_{z0}$. More precisely
${n_{z0}\over n_z}=\dot q_{z0}\pm O(n_z^{\smash{-1/2}})$ with probability
1 (w.p.1). Similarly for the right branch.
Assume the probabilities are equal ($\dot
q_{z0}=\dot q_{z1}=\odt$), possibly but not necessarily due to a
uniform $\dot q_z()$ on $\Gamma_z$. Then $\Delta_z=O(n_z^{\smash{-1/2}})$,
which implies
\beqn\label{eqwtoinfty}
  w_{n_z}(\Delta_z) \sim \Theta(\sqrt{n}_z)
  \;\mathop{\longrightarrow}\limits_{w.p.1}^{n_z\to\infty}\; \infty
  \qmbox{if} \dot q_{z0}=\dot q_{z1}=\odt,
\eeqn
consistent with our anticipation.
Conversely, for $\dot q_{z0}\neq\dot q_{z1}$ (which implies
non-uniformity of $\dot q_z()$) we have $\Delta_z\to c:=\dot
q_{z0}-\odt \neq 0$, which implies
\beqn\label{eqwto0}
  w_{n_z}(\Delta_z) \sim {\textstyle\sqrt{2n_z\over\pi}}\,\e^{-2n_z c^2}
  \;\mathop{\longrightarrow}\limits_{w.p.1}^{n_z\to\infty}\; 0
  \qmbox{if} \dot q_{z0}\neq\dot q_{z1},
\eeqn
again, consistent with our anticipation.

\paradot{\boldmath Asymptotic convergence/consistency ($n\to\infty$)}
For fixed $m<\infty$, one can show that almost surely the
posterior $p_z(\v q_{z*}|D)$ concentrates around the true
distribution $\,\dot{\!\v q}_{z*}$ for $n\to\infty$. This implies
that the posterior $p_z(x|D_z)\to\dot q_z(x)$ for all
$x\in\Gamma_z$. One can also show that the evidence $p_z(D_z)\to $
$\odt$ or $1$ for uniform $\dot q_z()$, and increases
exponentially with $n_z$ for non-uniform $\dot q_z()$ (see
\cite{Hutter:04btcode} for proofs).

\paradot{Model dimension and cell number}
As discussed in Section \ref{secTMM}, the effective dimension of
$\v q_*$ is the number of components that are not forced to $\odt$
by (\ref{tmmu1}). Note that a component may be ``accidentally''
$\odt$ in (\ref{tmms1}), but since this is an event of probability
0, we don't have to care about this subtlety. So the effective
dimension $N_{\v q_{z*}}=\#\{q\in\v q_{z*}:q\neq\odt\}$ of $\v
q_{z*}$ can be given recursively as
\beq\label{eqNrec}
   N_{\v q_{z*}} \;=\;
   \left\{ {0 \qquad\;\;\qmbox{if} \l(z)=m \qmbox{or} q_{z0}=\odt \atop
            1+N_{\v q_{z0*}}+N_{\v q_{z1*}} \;\;\qquad\qquad\qmbox{else}}
   \right.
\eeq
The effective dimension is zero if $q_{z}=\odt$, since this
implies that the whole tree $\Gamma_z$ has $q_{zy}=\odt$ due to
(\ref{EqPqRec}). If $q_z\neq\odt$, we add the effective dimensions
of subtrees $\Gamma_{z0}$ and $\Gamma_{z1}$ to the root degree of
freedom $q_{z0}=q_z-q_{z1}$.
Bayes' rule allows to represent the posterior probability that
$N_{\v q_{z*}}=k$ as
\beqn
  P_z[N_{\v q_{z*}}=k|D_z]\!\cdot\!p_z(D_z)
  = \!\int\! \delta_{N_{\v q_{z*}}k} p_z(D_z|\v q_{z*})p(\v q_{z*})d\v q_{z*}
\eeqn
where $P_z[...|...]:=P[...|\Gamma_z...]$, and
$\delta_{ab}=1$ for $a=b$ and $0$ else.
The r.h.s.\ coincides with
(\ref{eqEvRecProof}) except for the extra factor $\delta_{N_{\v q_{z*}}k}$.
Analogous to the evidence
(\ref{eqEvRecProof}), using (\ref{eqNrec}) we
can prove the following recursion:
\bqa\nonumber
  & & \nq\nq P_z[N_{\v q_{z*}}=0|D_z] \;=\;
   1-g_z(D_z), \qquad \qmbox{for} l<m,
\\[1ex] \label{eqMDR} \label{eqEMD}
  & & \nq\nq P_z[N_{\v q_{z*}}=k+1|D_z] \;=\;
\\[-1ex] \nonumber
  & & \nq\nq g_z(D_z)\!\cdot\!\!\sum_{i=0}^k
  P_{z0}[N_{\v q_{z0*}}\!=\!i|D_{z0}]\cdot P_{z1}[N_{\v q_{z1*}}\!=\!k\!-\!i|D_{z1}],
\\ \nonumber
  & & \nq\nq P_z[N_{\v q_{z*}}=k|D_z]=\delta_{k0}
  :=\textstyle\big\{{1 \;{\rm if}\; k=0 \atop 0 \;{\rm if}\; k>0}\big\} \qmbox{for} l=m.
\\ \label{eqSplitProb}
   & & \nq\nq g_z(D_z) := \odt{p_{z0}(D_{z0})p_{z1}(D_{z1})
                   \over p_z(D_z) w({n_{z0},n_{z1}}) }
  \stackrel{(\ref{eqEvDens})}= 1-{1\over 2 p_z(D_z)}
\eqa
Read: The probability that tree $\Gamma_z$ has dimension $k+1$
equals the posterior probability $g_z(D_z)$ of splitting
$\Gamma_z$, times the probability that left subtree has dimension
$i$, times the probability that right subtree has dimension $k-i$,
summed over all possible $i$.

Let us define a cell or bin as a maximal volume on which $q()$ is
constant. Then the model dimension is 1 less than the number of
bins (due to the probability constraint). Hence we also have a
recursion for the distribution of the number of cells.

\paradot{Tree height and cell size}
The effective height of tree $\v q_{z*}$ at $x\in\Gamma_z$
is also an interesting property.
If $q_{z0}=\odt$ or $\l(z)=m$, then the height $h_{\v q_{z*}}(x)$ of
tree $\v q_{z*}$ at $x$ is obviously zero. If
$q_{z0}\neq\odt$, we take the height of the subtree
$\v q_{zx_{l+1}*}$ that contains $x$ and add 1:
\beqn\label{eqhrec}
   h_{\v q_{z*}}(x) \;=\;
   \left\{ {0 \atop 1+h_{\v q_{zx_{l+1}*}}(x)}
           {\nq\mbox{if}\quad \l(z)=m \qmbox{or} q_{z0}=\odt \atop
             \mbox{else}\hspace{10ex}}   \right.
\eeqn
One can show that the tree height at $x$ averaged over all trees $\v q_{z*}$ is
\beqn
  E_z\![h_{\v q_{z*}}\!(x)|D_z] =
  g_z(D_z)\Big[1 + E_{zx_{l+1}}\![h_{\v q_{zx_{l+1}\!*\!}}(x)|D_{zx_{l+1}}]\Big]
\eeqn
where $E_z[f_{\v q_{z*}}|...]=\int P_z[f_{\v q_{z*}}|...]p(\v
q_{z*})d\v q_{z*}$. We may also want to compute the tree height
averaged over all $x\in\Gamma_z$. For $\l(z)<m$ and
$q_{z0}\neq\odt$ we get
\beqn
  \bar h_{\v q_{z*}}
  \;:=\;
  \int h_{\v q_{z*}}(x)q(x|\Gamma_z)dx
  \;=\; 1+q_{z0}\!\cdot\!\bar h_{\v q_{z0*}}+q_{z1}\!\cdot\!\bar h_{\v q_{z1*}}
\eeqn\vspace{-2ex}
\bqan
  E_z[\bar h_{\v q_{z*}}|D_z] \;=\;
  g_z(D_z)\Big[1 \nq\;&+&\nq\; {n_{z0}\!+1\over n_z+2}E_{z0}[\bar h_{\v q_{z0*}}|D_{z0}]
\\
            \nq\;&+&\nq\; {n_{z1}\!+1\over n_z+2}E_{z1}[\bar h_{\v q_{z1*}}|D_{z1}]\Big]
\eqan
with obvious interpretation: The expected height of a subtree is
weighted by its relative importance, that is (an estimate of) its
probability. The recursion terminates with $E_z[h_{\v
q_{z*}}|D_z]=0$ when $\l(z)=m$. We can also compute intra and
inter tree height variances.

Finally consider the average cell size or volume $v$. Maybe more
useful is to consider the logarithm $-\lb|\Gamma_z|=\l(z)$, since
otherwise small volumes can get swamped in the expectation by a
single large one. Log-volume $v_{\v q_{z*}}=\l(z)$ if $\l(z)=m$ or
$q_z=\odt$, and else recursively $v_{\v q_{z*}}=q_{z0}v_{\v
q_{z0*}}+q_{z1}v_{\v q_{z1*}}$. We can reduce this to the tree height,
since $v_{\v q_{z*}}=\bar h_{\v q_{z*}}+\l(z)$, in particular
$v_{\v q_*}=\bar h_{\v q_*}$

\section{\boldmath INFINITE TREES ($m\to\infty$)}\label{secIT}

\paradot{Motivation}
We have chosen an (arbitrary) finite tree height $m$ in our setup,
needed to have a well-defined recursion start at the leaves of the
trees. What we are really interested in are infinite trees
($m=\infty$).
Why not feel lucky with finite $m$? First, for continuous domain
$\Gamma$ (e.g.\ interval $[0,1)$), our tree model contains only
piecewise constant models. The true distribution $\dot q()$ is
typically non-constant and continuous (Beta, normal, ...). Such
distributions are outside a finite tree model class (but inside
the infinite model), and the posterior $p(x|D)$ cannot converge to
the true distribution, since it is also piecewise constant. Hence
all other estimators based on the posterior are also not
consistent. Second, a finite $m$ violates scale invariance (a
non-informative prior on $\Gamma_z$ should be the same for all $z$,
apart from scaling). Finally, having to choose the ``right'' $m$
may be worrisome.

For increasing $m$, the cells $\Gamma_x$ become smaller and will
(normally) eventually contain either only a single data item, or
be empty. It should not matter whether we further subdivide empty
or singleton cells. So we expect inferences to be independent of
$m$ for sufficiently large $m$, or at least the limit $m\to\infty$
to exist. In this section we show that this is essentially true.

\paradot{\boldmath Prior inferences ($D=\es$)}
We first consider the prior (zero data) case $D=\es$. Recall that
$z\in\SetB_0^m$ is some node and $x\in\SetB^m$ a leaf node.
Normalization implies $p_z(\es)=1$ for all $z$, which is
independent of $m$, hence the prior evidence exists for
$m\to\infty$. This is nice, but hardly surprising.

The prior effective model dimension $N_{\v q_*}$ is more
interesting. $D=\es$ implies $D_z=\es$ implies $n_z=0$ implies
$w(n_{z0},n_{z1})=1$ implies a 50/50 prior chance $g_z(\es)=\odt$
for a split (see (\ref{eqSplitProb})). Recursion (\ref{eqMDR})
reads
\beqn
  P_z[N_{\v q_{z*}}=k+1] =
  \odt\sum_{i=0}^k P_{z0}[N_{\v q_{z0*}}\!=\!i]\cdot P_{z1}[N_{\v q_{z1*}}\!=\!k\!-\!i]
\eeqn
with $P_z[N_{\v q_{z*}}=k]=\delta_{k0}$ for $l=m$ and $P_z[N_{\v
q_{z*}}=0]=\odt$ for $l<m$. So the recursion terminates in recursion depth
$\min\{k+1,m-l\}$. Hence $P_z[N_{\v q_{z*}}=k+1]$ is the same for
all $m>l+k$, which implies that the limit $m\to\infty$ exists.
Furthermore, recursion and termination are independent of $z$,
hence also $a_k:=P_z[N_{\v q_{z*}}=k]$. So we have to solve the recursion
\beq\label{eqMDar}
   a_{k+1} = \odt\sum_{i=0}^k a_i\!\cdot\! a_{k-i} \qmbox{with} a_0=\odt
\eeq
The first few coefficients can be bootstrapped by hand:
($\odt$,$1\over 8$,$1\over 16$,$5\over 128$,$7\over 256$,$21\over
1024$,$33\over 2048$,...).
A closed form  can
also be obtained: Inserting (\ref{eqMDar}) into $f(x):=\sum_{k=0}^\infty
a_k x^{k+1}$ we get $f(x)=\odt[x+f^2(x)]$ with solution
$f(x)=1-\sqrt{1-x}$, which has Taylor expansion coefficients
\beqn
   a_k = (-)^k\left({1/2\atop k+1}\right)
   = {1\over 2(k\!+\!1)4^k}\left({2k\atop k}\right) \sim {1\over 2\sqrt{\pi}}\; k^{-3/2}
\eeqn
$(a_k)_{k\in\SetN_0}$ is a well-behaved distribution. It decreases
fast enough to be a proper measure ($\sum_k a_k=f(1)=1<\infty$),
but too slow for the expectation $E[N_{\v q_*}]=\sum_k k\cdot
a_k=\infty$ to exist. This is exactly how a proper non-informative
prior on $\SetN$ should look like: as uniform as possible, i.e.\
slowly decreasing. Further, $P[N_{\v q_*}<\infty]=\sum_k a_k=1$
implies $P[N_{\v q_*}<\infty|D]=1$, which shows that the effective
dimension is almost surely finite, i.e. infinite (Polya) trees
have probability zero.

For the tree height we have $E_z[h_{\v q_{z*}}(x)]=0$ if $l=m$ and
otherwise
\bqan
  E_z[h_{\v q_{z*}}(x)] &=&
  \odt [1+E_{zx_{l+1}}[h_{\v q_{zx_{l+1}*}}(x)]]
\\
  &=&...\;= 1-(\odt)^{m-l} \to 1 \qmbox{for} m\to\infty
\eqan
This also implies that the expected average height $E_z[\bar h_{\v
q_{z*}}]=1-(\odt)^{m-l} \to 1$. This is the first case where the
result is not independent of $m$ for large finite $m$, but it
converges for $m\to\infty$, what is enough for our purpose.

\paradot{\boldmath Single data item $D=(x)$}
Since $p(x)\equiv 1$ (by symmetry and normalization) and $w_1=1$
are the same as for the $n=0$ case, all prior $n=0$, $m\to\infty$
results remain valid for $n=1$: $g(x)=\odt$, $P[N_{\v
q_*}=k|x]=a_k$, and $E[h_{\v q_*}(x)|x]\to 1$.

\paradot{\boldmath General $D$}
We now consider general $D$. For continuous spaces $\Gamma$ and
non-singular distribution $\dot q$, the probability of observing
the same point more than once (multi-points) is zero and hence
can, to a certain extend, be ignored. See \cite{Hutter:04btcode}
for a thorough workout of this case.
In order to compute $p(D)$ and other quantities, we recurse
(\ref{eqEvDens}) down the tree until $D_z$ is either empty or
a singleton $D_z=(x)\in\Gamma_z$. We call the depth $m_x:=\l(z)$
at which this happens, the separation level. In this way, the
recursion always terminates. For instance, for $\Gamma=[0,1)$, if
$\eps:=\min\{|x^i-x^j|:x^i\neq x^j$ with $x^i,x^j\in D\}$ is the
shortest distance, then $m_x<\lb{2\over\eps}=:m_0<\infty$,
since $\eps>0$. At the separation level we can insert the derived
formulas for evidence, posterior, dimension, and height.
Note, there is no approximation here. The procedure
is exact, since we analytically computed the infinite recursion
for empty and singleton $D$.

So we have devised a finite procedure, linear in the data size
$n$, for exactly computing all quantities of interest in the
infinite Bayes tree. In the worst case, we have to recurse down to
level $m_0$ for each data point, hence our procedure has
computational complexity $O(n\cdot m_0)$. For non-singular prior,
the time is actually $O(n)$ with probability 1. So, inference in
our mixture tree model is {\em very} fast.
Posterior (weak) convergence/consistency for $m=\infty$ can be
shown similarly to the $m<\infty$ case \cite{Hutter:04btcode}.

\section{THE ALGORITHM}\label{secAlg}

\paradot{What it computes}
In the last two sections we derived all necessary formulas for
making inferences with our tree model. Collecting pieces together
we get the exact algorithm for infinite tree mixtures below. It
computes the evidence $p(D)$, the expected tree height $E[h_{\v
q_*}(x)|D]$ at $x$, the average expected tree height $E[\bar h_{\v
q_*}|D]$, and the model dimension distribution $P[N_{\v q_*}|D]$.
It also returns the number of recursive function calls, i.e.\ the
size of the explicitly generated tree. The size is proportional to
$n$ for regular distributions $\dot q$.

\paranodot{The BayesTree algorithm} (in pseudo C code)
takes arguments $(D[],n,x,N)$; data array $D[0..n-1]\in[0,1)^n$, a
point $x\in\SetR$, and an integer $N$. It returns $(p,h,\bar
h,\tilde p[],r)$; the logarithmic data evidence $p\widehat=\ln p(D)$, the
expected tree height $h\widehat=E[h_{\v q_*}(x)|D]$ at $x$, the average expected tree
height $\bar h\widehat=E[\bar h_{\v q_*}|D]$, the model dimension distribution $\tilde
p[0..N-1]\widehat=P[N_{\v q_*}=..|D]$, and the number of recursive function
calls $r$ i.e.\ the size of the generated tree. Computation time is
about $N^2 n\log n$ nano-seconds on a 1GHz P4 laptop.

\begin{list}{}{\parskip=0ex\parsep=0ex\itemsep=0.5ex\leftmargin=0ex\labelwidth=0ex}
  \item {\bf\boldmath BayesTree($D[],n,x,N$)}
  \begin{list}{}{\parskip=0ex\parsep=0ex\itemsep=0.5ex\leftmargin=2ex\labelwidth=1ex\labelsep=1ex}
    \item[$\lceil$] if ($n\leq 1$ and ($n==0$ or $D[0]==x$ or $x\not\in[0,1)$))
    \begin{list}{}{\parskip=0ex\parsep=0ex\itemsep=0.5ex\leftmargin=2ex\labelwidth=1ex\labelsep=1ex}
      \item[$\lceil$] if ($x\in[0,1)$) then $h=1$; else $h=0$;
      \item     $\bar h=1$; $p=\ln(1)$; $r=1$;
      \item[$\lfloor$] for$(k=0,..,N-1)$ $\tilde p[k]=a_k$; \hfill /* see (\ref{eqMDar}) */
    \end{list}
      \item else
    \begin{list}{}{\parskip=0ex\parsep=0ex\itemsep=0.5ex\leftmargin=2ex\labelwidth=1ex\labelsep=1ex}
      \item[$\lceil$] $n_0=n_1=0$;
      \item           for$(i=0,..,n-1)$
      \begin{list}{}{\parskip=0ex\parsep=0ex\itemsep=0.5ex\leftmargin=2ex\labelwidth=1ex\labelsep=1ex}
        \item[$\lceil$] if ($D[i]<\odt$) then[$\,D_0[n_0]=2D[i]$; \ \ \ $\,n_0=n_0+1$;]
        \item[$\lfloor$] \hspace{12ex} else [$D_1[n_1]=2D[i]-1$; $n_1=n_1+1$;]
      \end{list}
      \item ($p_0,h_0,\bar h_0,\tilde p_0[],r_0$)=BayesTree($D_0[],n_0,2x,N-1$);
      \item ($p_1,h_1,\bar h_1,\tilde p_1[],r_1$)=BayesTree($D_1[],n_1,2x\!-\!1,N\!-\!1$);
      \item $t=p_0+p_1-\ln w(n_0,n_1)$;
      \item if ($t<100$) then $p=\ln(\odt(1+\exp(t))$;
      \item \hspace{10ex} \ else \ $p=t-\ln(2)$;
      \item $g=1-\odt\exp(-p)$;
      \item if ($x\in[0,1)$) then $h=g\cdot(1+h_0+h_1)$; else $h=0$;
      \item $\bar h=g\cdot(1+{n_0+1\over n+2}\bar h_0+{n_1+1\over n+2}\bar h_1)$;
      \item $\tilde p[0]=1-g$;
      \item for($k=0,..,N-1$) $\tilde p[k+1]=g\cdot\!\sum_{i=0}^k \tilde p_0[i]\cdot \tilde p_1[k-i]$;
      \item[$\lfloor$] $r=1+r_0+r_1$;
    \end{list}
    \item [$\lfloor$] {\bf\boldmath return ($p,h,\bar h,\tilde p[],r$); }
  \end{list}
\end{list}

\paradot{How algorithm BayesTree() works}
Since evidence $p(D)$ and weight $1/w_n$ can grow exponentially
with $n$, we have to store and use their logarithms.
So the algorithm returns $p\widehat=\ln p(D)$.
In the $n\leq 1$ branch, the closed form solutions $p\widehat=\ln
p(\es)=\ln(1)$, $h\widehat=E[h_{\v q_*}(x)|\es\mbox{ or }x]=1$,
$\bar h\widehat=E[\bar h_{\v q_*}|D]=1$, and $\tilde
p[k]=a_k$ have been used to truncate the recursion. If $D=(x^1)\neq
x$, we have to recurse further until $x$ falls in an empty
interval.
In this case or if $n>1$ we partition $D$ into points left and
right of $\odt$. Then we rescale the points to [0,1) and store
them in $D_0$ and $D_1$, respectively. Array $D$ could have been
reused (like in quick sort) without allocating two new arrays.
Then, algorithm BayesTree() is recursively called for each
partition. The results are combined according to the recursions
derived in Section \ref{secTMM}.
$\ln w$ can be computed from (\ref{eqWeights}) via $\ln
n!=\sum_{k=1}^n\ln k$. (Practically, pre-tabulating $a_k$ or $n!$
does not improve overall performance).
For computing $p$ we need to use $\ln(\odt(1+e^t))\dot=t-\ln 2$ to
machine precision for large $t$ in order to avoid numerical overflow.

\paradot{Remarks}
Strictly speaking, the algorithm has runtime $O(n\log n)$, since
the sorting effectively runs once through all data at each level.
If we assume that the data are presorted or the counts $n_z$ are
given, then the algorithm is $O(n)$ \cite{Hutter:04btcode}. The
complete C code, available from \cite{Hutter:04btcode}, also
handles multi-points.

Note that $x$ passed to BayesTree() is {\em not} and cannot be
used to compute $p(x|D)$. For this, one has to call BayesTree()
twice, with $D$ and $(D,x)$, respectively. The quadratic order in
$N$ is due to the convolution, which could be reduced to $O(N\log
N)$ by transforming it to a scalar product in Fourier space with
FFT.

Multiply calling BayesTree(), e.g.\ for computing the predictive
density function $p(x|D)$ on a fine $x$-grid, is inefficient. But
it is easy to see that if we once pre-compute the evidence
$p_z(D_z)$ for all $z$ up to the separation level in time $O(n)$,
we can compute ``local'' quantities like $p(x|D)$ at $x$ in time
$O(\log n)$. This is because only the branch containing $x$ needs
to be recursed, the other branch is immediately available, since it
involves the already pre-computed evidence only. The predictive
density $p(x|D)=E[q(x)|D]$ and higher moments, the distribution
function $P[x\leq a|D]$, updating $D$ by adding or removing one
data item, and most other local quantities can be computed in time
$O(\log n)$ by such a linear recursion.

A good way of checking correctness of the implementation {\em and}
of the derived formulas, is to force some {\em minimal} recursion
depth $m'$. The results must be independent of $m'$, since the
closed-form speedups are exact and applicable anywhere beyond the
separation level.

\begin{figure*}
\centerline{
\includegraphics[width=0.33\textwidth]{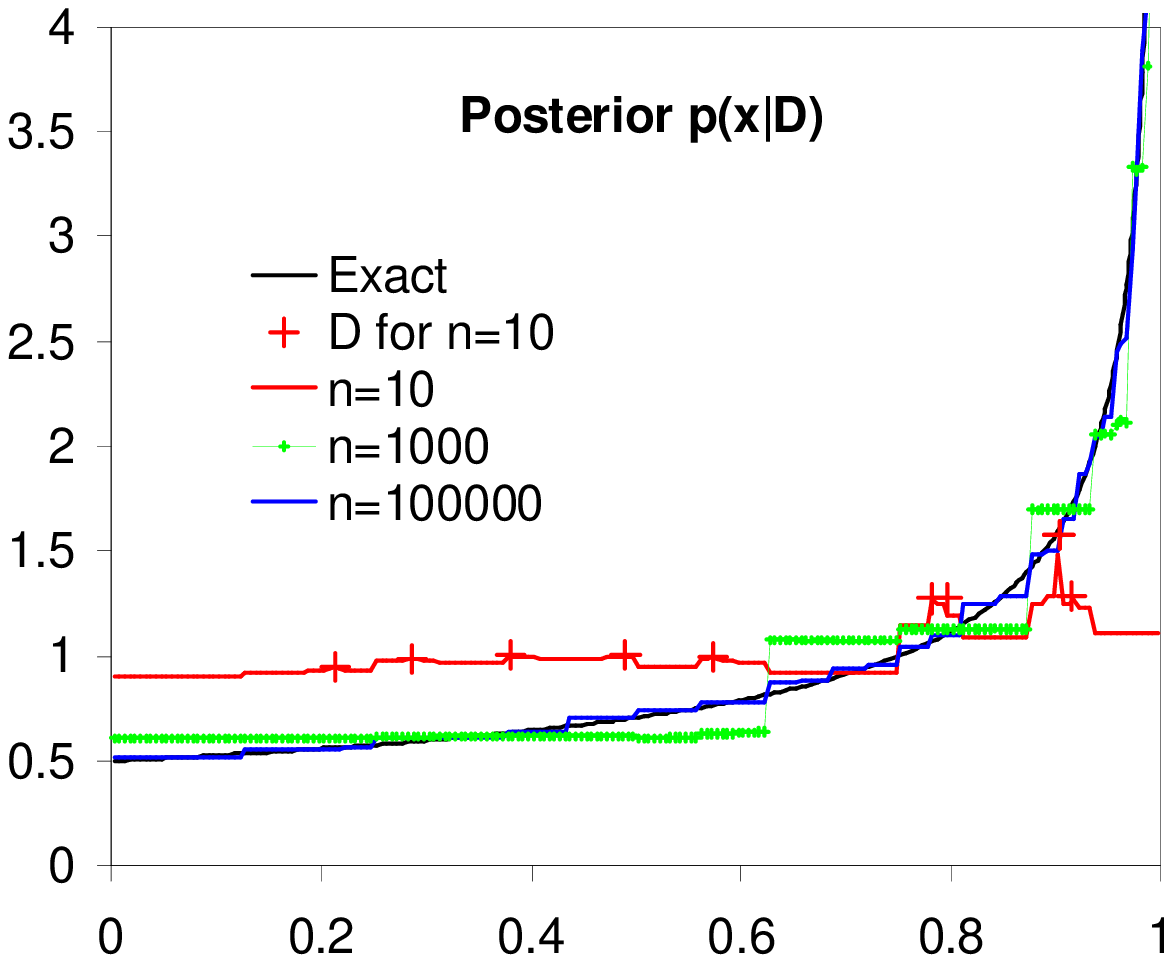}
\includegraphics[width=0.33\textwidth]{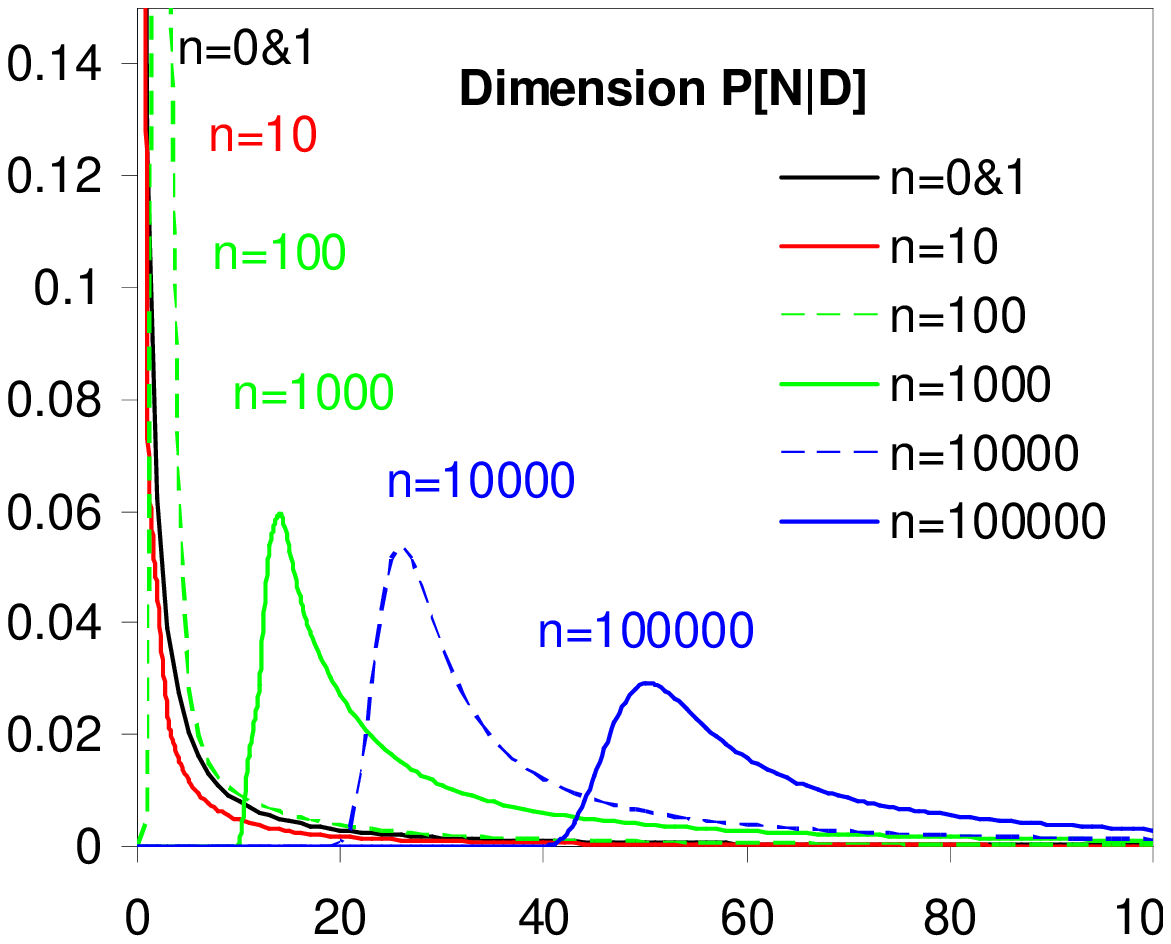}
\includegraphics[width=0.33\textwidth]{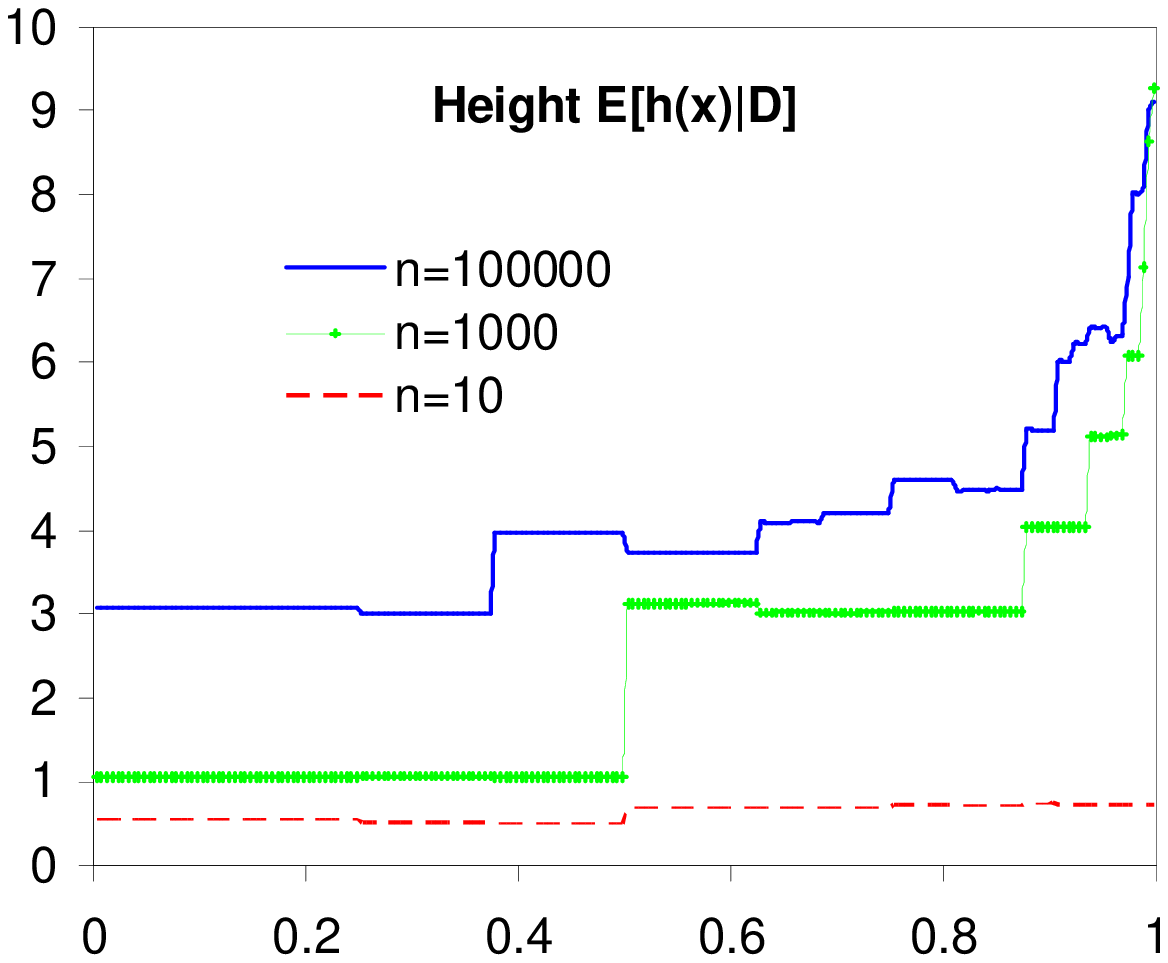} }
\caption{\label{figSingular1}BayesTree() results for a
prototypical proper singular distribution $\dot q(x)=2/\sqrt{1-x}$.}
\end{figure*}

\paradot{Numerical example}
To get further insight into the behavior of our model, we
numerically investigated some example distributions $\dot q()$. We
have chosen elementary functions, which can be regarded as
prototypes for more realistic functions. They include the Beta,
linear, a singular, piecewise constant distributions with finite
and infinite Bayes trees, and others.
These examples on $[0,1)$ also shed light on the other spaces
discussed in Section \ref{secTMM}, since they are isomorphic.
The posteriors, model dimensions, and tree heights, of
the singular distribution $\dot q(x)=2/\sqrt{1-x}$
are plotted in Figure \ref{figSingular1} for random
samples $D$ of sizes $n=10^0,...,10^5$.
The posterior $p(x|D)$ clearly converges for $n\to\infty$ to the
true distribution $\dot q()$, accompanied by a (necessary)
moderate growth of the effective dimension. For $n=10$ we show the
data points. It is visible how each data point pulls the posterior
up, as it should be (``one sample seldom comes alone'').
The expected tree height $E[h(x)|D]$ correctly reflects the local
needs for (non)splits, i.e.\ is larger near the singularity at
$x=1$.
The other examples display a similar behavior (see
\cite{Hutter:04btcode}).

\section{DISCUSSION}\label{secDisc}

We presented a Bayesian model on infinite trees, where we split a
node into two subtrees with prior probability $\odt$, and uniform
choice of the probability assigned to each subtree. We
devised closed form expressions for various inferential quantities
of interest at the data separation level, which led to an exact
algorithm with runtime essentially linear in the data size.
%
The theoretical and numerical model behavior was very reasonable,
e.g.\ consistency (no underfitting) and low finite effective dimension
(no overfitting).

There are various natural generalizations of our model.
The splitting probability $p(s)$ could be chosen different from $\odt$,
$k$-ary trees could be allowed, and the uniform prior over
subtrees could be generalized to Beta/Dirichlet distributions.
We were primarily interested in the case of zero prior knowledge,
hence zero model (hyper)parameters, but the generalizations above
make the model flexible enough, in case prior knowledge needs to
be incorporated.
The dependency on $p(s)$ is particularly interesting \cite{Hutter:04btcode}.
The expected entropy can also be computed by allowing fractional
counts $n_z$ and noting that $x\ln x = {d\over dx}
x^\alpha|_{\alpha=1}$ \cite{Hutter:01xentropy}.
A sort of maximum a posteriori (MAP) tree skeleton can also easily
be read off from (\ref{eqEvDens}). A node $\Gamma_z$ in the
MAP-like tree is a leaf iff ${p_{z0}(D_{z0})p_{z1}(D_{z1})\over
w(n_{z0},n_{z1})}<1$.
A challenge is to generalize the model from piecewise constant to
piecewise linear continuous functions, at least for
$\Gamma=[0,1)$. Independence of subtrees no longer holds, which
was key in our analysis.


\bibliographystyle{alpha}
\begin{small}

\end{small}

\end{document}